\documentclass[a4paper,12pt]{article}
\usepackage{amsmath,amsfonts,enumerate,amssymb,amsthm,color}

\title{The number of cyclic configurations of type $(v_3)$ and the isomorphism problem}

\author{
S. Hiroki  Koike-Quintanar
\footnote{Supported in part by ARRS - Agencija za raziskovanje Republike Slovenija, program no. P1-0285.} \\   
{\small IAM, University of Primorska} \\  [-0.8ex]
{\small  Muzejski trg 2, 6000 Koper, Slovenia} \\  [-0.8ex]
{\small \texttt{hiroki.koike@iam.upr.si}} \\ [+1.5ex]
Istv\'an Kov\'acs
\footnote{Supported in part by ARRS - Agencija za raziskovanje Republike Slovenija, program no. P1-0285 and project N1-0011 (ESF EUROCORES EUROGiga/GReGAS)..} \\   
{\small IAM and FAMNIT, University of Primorska}  \\  [-0.8ex]
{\small  Muzejski trg 2, 6000 Koper, Slovenia} \\  [-0.8ex]
{\small  \texttt{istvan.kovacs@upr.si}} \\ [+1.5ex]
Toma\v z  Pisanski
\footnote{Supported in part by ARRS - Agencija za raziskovanje Republike Slovenija, program no. 
P1-0294 and project N1-0011 (ESF EUROCORES EUROGiga/GReGAS).}  \\ 
{\small FMF, University of Ljubljana} \\ [-0.8ex]
{\small  Jadranska 19, 1000 Ljubljana, Slovenia} \\ [-0.8ex] 
{\small \texttt{Tomaz.Pisanski@fmf.uni-lj.si}}
}

\def\B{\mathcal{B}}
\def\C{\mathcal{C}}

\def\K{\mathcal{K}} 
\def\mo{\mathrm{ord_m}}

\def\N{\mathbb{N}}
\def\proof{\noindent{\sc Proof:}\ }
\def\QED{\hfill\vrule height .9ex width .8ex depth -.1ex \medskip}
\def\Z{\mathbb{Z}}
\DeclareMathOperator{\aut}{Aut}
\DeclareMathOperator{\con}{Con}
\DeclareMathOperator{\obj}{Obj}
\DeclareMathOperator{\orb}{Orb}

\newtheorem{thm}{Theorem}[section]
\newtheorem{lem}[thm]{Lemma}

\newtheorem{cor}[thm]{Corollary}
\newtheorem{que}[thm]{Question}
\newtheorem*{thmA}{Theorem A}
\newtheorem*{thmB}{Theorem B}

\theoremstyle{definition}

\theoremstyle{remark}
\newtheorem{rem}[thm]{Remark}

\begin{document}

\maketitle

\begin{abstract}
A configuration of points and lines is cyclic if it has an automorphism which permutes 
its points in a full cycle.  A closed formula is derived for the number of non-isomorphic connected cyclic 
configurations of type $(v_3),$ i.e., which have $v$ points and lines, and each point/line is incident 
with exactly  $3$ lines/points.  
In addition, a Bays-Lambossy type theorem is proved for cyclic configurations 
if the number of points is  a product of two primes or a prime power.

\medskip\noindent{\it Keywords:} cyclic configuration, cyclic object,  isomorphism.

\medskip\noindent{\it MSC 2010:}  20B25, 51E30, 05C25, 05C60.
\end{abstract}

\section{Introduction}

An {\em incidence geometry} $(P,\B)$ consists of a set of $v$ points $P = \{p_1,...,p_v\}$ and a collection of $b$ 
lines (or blocks) $\B = \{B_1,...,B_b\}$ such that $B_i \subseteq P$ for every $i \in \{1,...,b\},$ and 
$|B_i \cap B_j| \le 1$ for every $i,j \in \{1,...,b\}$ and $i \ne j$. An incidence geometry is called a {\em 
configuration} of type $(v_r,b_k)$ ({\em combinatorial configuration} in the sense of \cite{Gru09}) if 
\begin{itemize}
\item $| \{ B_j \in \B : p_i \in B_j \} | = r$ for every $i \in \{1,..,v\};$ and 
\item $|B_j| = k$ for every $j \in \{1,..,b\}$ with $k \ge 3$.
\end{itemize}
A configuration with $v = b$ (and therefore $r = k$) is called \textit{balanced,} 
or a \textit{$k$-configuration,}  
and its type is simply denoted by $(v_k)$. 
A configuration $(P,\B)$ is called \textit{decomposable} if it is the disjoint union of two 
configurations $(P_r,\B_r), r=1,2,$  i.e.,  $P = P_1 \cup P_2,$ $P_1 \cap P_2 = \emptyset,$ and 
$\B = \B_1 \cup \B_2$.   Indecomposable configurations are also called \textit{connected}. 
An \textit{isomorphism} between two incidence geometries  $(P_r,\B_r), r = 1,2,$ is a bijective mapping 
$\sigma : P_1 \to P_2$  which maps $\B_1$ onto $\B_2$. Here a block $B \in \B_1$ with $B = \{p_1,...,p_k\}$ 
is mapped onto $B^\sigma = \{p_1^\sigma,...,p_k^\sigma\}$.  If $(P_1,\B_1) = (P_2,\B_2),$ then 
$\sigma$ is called an \textit{automorphism}; the group of all automorphisms will be denoted by 
$\aut((P_1,\B_1))$.   
An incidence geometry is \textit{cyclic},  if it has an 
automorphism  which permutes its points  in a full cycle.  In this case it is natural to identify the points with  elements in the 
ring $\Z_v,$ and assume that the translation $x \mapsto x + 1$  is an automorphism. Now, two incidence geometries are said to be \textit{multiplier equivalent},  if there exits a unit $a \in \Z_n^*$ such that the 
mapping $x \mapsto a x$ induces an isomorphism between them. 

The enumeration problem for configurations (both geometrical and combinatorial) attracted considerable 
attention (see the monograph \cite[Chapters 2-3]{Gru09}). The list of all configurations of type $(v_3)$ up to $v = 18$ was 
produced in \cite{BetBP00}, and for an approach based on the respective Levi graphs, we refer to  
\cite{BobPZ05,PetP05,PetZ09}. 
The latter approach is based on the easy but crucial observation that combinatorial $k$-configurations are the same things as bipartite $k$-valent graphs with a given black-and-white coloring. 
In this paper we are going to calculate the number of cyclic configurations of type $(v_3)$. 
For this purpose we set the notation $\# C(v_k)$ for the number of non-isomorphic coonected cyclic configurations  of 
type $(v_k)$. 
Our main result is the following closed formula for $\#C(v_3)$:

\begin{thmA}
Let  $v > 4$ be an integer with prime factorization $v = p_1^{n_1} \cdots p_k^{n_k}$. Then 
\begin{equation}\label{formula}
 \# C(v_3) = \left\{ \begin{array}{ll}
\frac{v}{6} \prod_{i=1}^k \big(1+\frac{1}{p_i} \big)  +  \alpha 2^k - 2  & \mbox{ if $v$ is odd},  \\ [+1.5ex]
\frac{v}{6}  \prod_{i=1}^k \big(1+\frac{1}{p_i} \big) +  \beta 2^k  - 3   & \mbox{ if $v$ is even},
\end{array} \right. 
\end{equation}
where $\alpha$ is defined for $v$ odd by 
\[
\alpha = \begin{cases} 5/6  &  \text{if   every } p_i \equiv 1(\text{mod }3), \\ 
                                      2/3  &  \text{if } p_1^{n_1} = 3   
                                             \text{ and if $i > 1,$ then } p_i \equiv 1 (\text{mod }3), \\
                                      1/2  &  \text{otherwise}, 
\end{cases} 
\]
and $\beta$ is defined for $v$ even by
\[
\beta  = \begin{cases} 1/4  &  \text{if  \ } v \equiv 2 (\text{mod }8) \text{ or } v \equiv 6 (\text{mod }8), \\ 
                                     1/2  &  \text{if  \ } v \equiv 4 (\text{mod }8), \\
                                        1  &  \text{if  \ } v \equiv 0 (\text{mod }8).
\end{cases}
\]
\end{thmA}

The crucial fact towards Theorem A is that the  isomorphism  problem in this case can be solved 
entirely by means of multiplier equivalence. More precisely,  every two cyclic configurations of type $(v_3)$ are 
isomorphic if and only if they are multiplier equivalent. 
This fact we are going to deduce as a direct consequence 
of a  result  about circulant matrices proved in \cite{WieZ07}. 
It is interesting to note that this is no longer true for arbitrary cyclic configurations with $3$ points on 
a line. In \cite{Phe87}, Phelps gives examples of cyclic $2$--$(v,3,1)$  designs which are  isomorphic but 
not multiplier equivalent. 
In Section 2 we review the relevant results on circulant  matrices and explain the relation with 
configurations. Section 3 is devoted to the proof of Theorem A. 

\medskip

In Section 4 we turn to the following question: 

\begin{que}\label{Q} 
Given an integer $v,$ is it true that any two balanced cyclic configurations on $v$ points
are isomorphic if and only if they are multiplier equivalent? 
\end{que} 

This is part of the more general question which asks if a given finite group $G$ has the \textit{CI-property} for a given class 
$\K$ of relational structures on $G$ (see \cite{Bab77,Pal87}).  This question  has been extensively  studied  
under various choices of $G$ and $\K$ (see, e.g. \cite{Huf96,HufJP93,Jun08,Muz99,Phe87}, 
just to mention those papers that will be invoked in the sequel). 
The finite groups having the CI-property for all relational structures (for short the CI-groups) were characterized by 
P\'alfy in \cite{Pal87}.  It turns out that these are the cyclic groups of order $n$ with $n=4$ or $\gcd(n,\phi(n)) = 1$ where $\phi$ denotes Euler's $\phi$  function. Consequently, Question \ref{Q} is answered in the positive 
if $v \ge 7$ and $\gcd(v,\phi(v)) = 1$. In Section 4 we provide further values of $v$  inducing a 
positive answer by proving the following theorem: 

\begin{thmB}
If  $v = pq$  or $v = p^n$, $p,q$ are primes, then any two balanced cyclic configurations on $v$ points 
are  isomorphic if and only if they are multiplier equivalent.
\end{thmB}

\begin{rem}
Theorem B can be viewed as a Bays-Lambossy type theorem for balanced cyclic configurations. 
It was proved first by Bays in \cite{Bay30} and Lambossy in \cite{Lam31} that two cyclic Steiner triple systems 
on a prime number of points are isomorphic if and only if they are multiplier equivalent. 
It is worth noted that the Bays-Lambossy Theorem was generalized to abelian projective planes (see \cite[Corollary 2.2]{Jun08}).
\end{rem}

\section{Circulant matrices}

\begin{lem}
Let $\C = (\Z_v,\B)$  be a balanced configuration with the translation $x \mapsto x+1$  in  
$\aut(\C)$.   Then there exists a subset $S$ of $\Z_v$ such that $\B$ consists of the sets in the form  $S + i,$ $i \in \Z_v$. 
\end{lem}

\proof Denote by $G$ the group generated by the translation $x \mapsto x+1$. 
Choose a line $B \in \B$ such that $0 \in B,$ where $0$ is the zero element of $\Z_v$. Assume for the moment that $B$ satisfies 
\begin{equation}\label{block}
B^g = B \text{ or } B^g \cap B = \emptyset \text{ for every } g \in G.
\end{equation}

In other words, $B$ is a block for the permutation group $G$ (see \cite[page 12.]{DixM96}). 
This gives that $B$ is an orbit of a subgroup of $G$ of size $k$ (see \cite[Theorem 1.5A]{DixM96}), where $k$ is the size of the lines.   
Since $G$ is a cyclic group, the set $B$ is uniquely determined. Choose next a line 
$B' \in \B$ for which  $0 \in B'$ and $B' \ne B$.  Then \eqref{block}  does not hold for $B',$ i.e., there exists $g \in G$ such that 
$B'$ and $B'^{\; g}$ intersects at a unique point, say $i$ ($i \in \Z_v$). 

Let us consider the action of $G$ on the set $\B$. We denote by $G_{B'}$ the stabilizer of the line $B'$ in this action, i.e., 
$G_{B'} = \{g \in G : B'^{\; g} = B'\}$. 
Then $G_{B'^{\; g}} = g^{-1} G_{B'} g,$ $|G_{B'^{\; g}}| = |G_{B'}|,$ and so $G_{B'^{\; g}} = G_{B'}$ (again, 
$G$ is a cyclic group). Clearly, every element in $G_{B'} \cap G_{B'^{\; g}}$ fixes the point $i$. 
Since $G$ is regular on the points, we obtain $G_{B'} = G_{B'}  \cap G_{B'^{\; g}} = 1$. 
The orbit-stabilizer property  (see \cite[Theorem 1.4A]{DixM96}) gives that 
the orbit of $B'$ under $G$  is of length $|G| = |P| = |\B|$. Letting $S = B',$ the lemma follows. 
\QED

We shall refer to the set $S$ in  Lemma 6 as a \textit{base line} of  $\C,$ and use the symbol $\con(\Z_v,S)$ for 
$\C$.  Base lines are characterized in the next lemma.

\begin{lem}\label{HMP} 
{\rm \cite{HlaMP02}} The following (1)-(2) are equivalent for every subset $S$ of $\Z_v$.
\begin{enumerate}[(1)]
\item $S$  is a  base line of a cyclic configuration of type $(v_k)$.
\item $|S| = k$ and $|S - S| = k^2-k+1$.\footnote{Here $S-S = \{ s_1 - s_2 : s_1,s_2 \in S\}$.}  
\end{enumerate}
\end{lem}

Suppose that $S$ is a base line such that $0 \in S$ (clearly, every configuration admits base lines with this property). 
The set $S$ generates a subgroup of $\Z_v,$ say of order $d,$ and denote it by $\Z_d$. 
Then $\con(\Z_d,S)$  is a connected configuration. Also,  $\con(\Z_v,S)$ 
decomposes to the union of $v/d$ copies of $\con(\Z_d,S)$:
 \begin{equation}\label{decomp}
\con(\Z_v,S) \cong \con(\Z_d,S) \cup \cdots \cup \con(\Z_d,S). 
\end{equation}
Note that, if $S$ is an arbitrary base line ($0$ is not necessarily in $S$), then it holds:
\begin{equation}\label{conn}
\con(\Z_v,S) \text{ is connected } \iff \langle S-S \rangle = \Z_v.
\end{equation}  

The following necessary condition for a set to be a base line will be used frequently through the paper.  
It follows promptly from the second part in (2) of Lemma \ref{HMP}.

\begin{cor}\label{corHMP}
If a subset $S$ of $\Z_v$ is a base line of a cyclic configuration, then $S$ contains no $H$-coset 
for every nontrivial subgroup $H \le \Z_v$.  
\end{cor}

For positive integers $v$ and $k$ denote by $B(v,k)$ the set of all base lines of $\Z_v$ of 
size $k$, and by $B_{\rm con}(v,k)$ the set of those which define connected configurations. More  
formally,   
\[
\begin{array}{lcl}
B(v,k)     &=& \big\{  X \subseteq \Z_v  :  \; |X| = k  \text{ and } |X - X| = k^2 -  k + 1  \big\}, \\ [1.5ex]
B_{\rm con}(v,k) &=& \big\{  X \in B(v,k)  :  \langle X- X \rangle = \Z_v \big\}.
\end{array}
\]
Notice that, if $X \in B(v,k),$ $a \in \Z_v^*$ and $b \in \Z_v,$ then 
the set $a X+ b$ is also in $B(v,k)$. Hence the mapping $X \mapsto aX +b$ defines 
an action of  the \textit{affine group} $AGL_1(v)$  on $B(v,k)$. Clearly, the subset $B_{\rm con}(v,k)$ of $B(v,k)$  is invariant with respect to this action.

\medskip

Next, we review the definition of a circulant matrix. 
Let $A$ be an $v$-by-$v$ matrix. The matrix $A$ is a \textit{permutation matrix} if it is a $(0,1)$ matrix, and  
every row and column contains exactly one $1$'s. Furthermore, $A = (a_{i, j})$ is a \textit{circulant matrix} if 
$a_{i+1,j+1} = a_{i, j}$ holds for every $i,j \in \{0,1,...,v-1\},$ where the additions in subscripts are modulo $v$. 
Here we label rows and columns by elements of $\Z_v$. We let $\Z_v = \{0,1,...,v-1\},$ the  leftmost column is labeled 
$0,$ the next is $1$ and so on. 
If $A = (a_{i,j})$ is an $v$-by-$v$ $(0,1)$  circulant matrix, then denote by $S_A$ the subset of $\Z_v$ defined by 
\[
S_A = \{ i \in  \Z_v : a_{0, \, i} =1 \}.
\]
The cardinality $|S_A|$ is also called the \textit{weight} of $A$. 
Also, $A^T$ denotes the transpose of the matrix $A$. 

Let $S \in B(v,k),$ and let $A$ be the $(0,1)$ circulant matrix defined by $S_A = S$. 
Then it follows immediately from the definitions that, $A$ is a \textit{line-point incidence matrix} of the cyclic configuration 
$\con(\Z_v,S)$ (see \cite{Gru09}). 

\begin{lem}\label{LEM-matrix}
For $r = 1,2,$ let $S_r \in B(v,k),$ and let $A_r$ be the $(0,1)$ circulant matrix defined by $S_{A_r} = S_r$. 
The following equivalence holds: 
\[ 
\con(\Z_v,S_1) \cong \con(\Z_v,S_2) \iff A_1 = P A_2 Q 
\]
for some $v$-by-$v$ permutation matrices $P$ and $Q$. 
\end{lem}

\proof Let $P$ and $Q$ arbitrary $v$-by-$v$ permutation matrices. 
Associate then the permutation $\pi$ of $\Z_v$ with $P$ and the  permutation $\sigma$ of $\Z_v$ with $Q$ as follows:
\[
i^\pi = j \stackrel{def}{\iff} P_{i,j} = 1 \text{ and } i^\sigma = j \stackrel{def}{\iff} Q_{j,i} = 1 \text{ for evey } i,j \in \Z_v.
\]

Then 
\[
(PA_2Q)_{i,j} = \sum_{k,l  = 0}^{v-1} P_{i,k} (A_2)_{k,l} Q_{l,j} = (A_2)_{i^\pi, j^\sigma}. 
\]

Now, $A_1 = PA_2Q$ can be interpreted as the permutation $\sigma$ maps the line
$S_1+i$ to the line $S_2+i^\pi$. Equivalently, $\sigma$ induces an isomorphism from  
$\con(\Z_v,S_1)$ to $\con(\Z_v,S_2)$.  The lemma follows. 
\QED

Lemma \ref{LEM-matrix}  brings us to the following result of Wiedman and Zieve:

\begin{thm}\label{WZ} 
{\rm \cite[Theorem 1.1]{WieZ07}} The following (1)-(4) are equivalent for every two $v$-by-$v$ $(0,1)$ circulant matrices 
$A_1$ and $A_2$ of weight at most $3$.
\begin{enumerate}[(1)]
\item There is $a \in \Z_v^* $ and $b \in \Z_v$ such that $S_{A_1} = a S_{A_2} + b$.  
\item There are $v$-by-$v$ permutation matrices $P, Q$ such that $A_1 = P A_2 Q$.
\item There is an $v$-by-$v$ permutation matrix $P$ such that $A_1A_1^T = P A_2A_2^T P^{-1}$. 
\item The complex matrices $A_1 A_1^T$ and $A_2 A_2^T$ are similar. 
\end{enumerate}
\end{thm}

Notice that, the configurations $\con(\Z_v,S_1)$ and $\con(\Z_v,S_2)$ are multiplier equivalent if and only if  
$S_1 = aS_2+b$ for some $a \in \Z_v^*$ and $b \in \Z_v$. Combining this with Lemma \ref{LEM-matrix} and Theorem \ref{WZ}, 
we obtain the required equivalence for configurations of type $(v_3)$:

\begin{cor}\label{corWZ}
Any two cyclic configurations of type $(v_3)$ are isomorphic if and only if these are multiplier equivalent. 
\end{cor}

As pointed out in \cite{WieZ07}, the equivalences in Theorem \ref{WZ} do not hold when the weight $k \ge 4$. 
The following theorem settles the case $k = 4$. It is was  proved by the first two authors
in the context of cyclic Haar graphs (see \cite[Theorem 1.1]{KoiK}), below it is rephrased in terms of 
circulant matrices.

\begin{thm}\label{KK} 
The following (1)-(2) are equivalent for two $v$-by-$v$ $(0,1)$ circulant matrices $A_1$ and $A_2$ of weight $4$ such that 
$\langle S_{A_r} - S_{A_r} \rangle = \Z_v$ for both $r = 1,2$. 
\begin{enumerate}[(1)]
\item There exist $a_1,a_2 \in \Z_v^*$ and $b_1,b_2 \in \Z_v$ such that 
\begin{enumerate}[({1}a)]
\item $S_{A_1} = a_1 S_{A_2} + b_1 ;$  or 
\item $a_1 S_{A_1} + b_1 = \{ 0,x,y,y+u \}$  and $a_2 S_{A_2} + b_2 = \{ 0,x+u,y,y+u\},$ where $v=2u,$ $
\Z_v = \langle x,y \rangle,$ $2 \mid x,$ $(2x) \mid u$ and $x/2 \not\equiv y +  u/(2x) (\text{mod } u/x)$.
\end{enumerate}
\item There are $v$-by-$v$ permutation matrices $P,Q$ such that $A_1 = P A_2 Q$. 
\end{enumerate}
\end{thm}

\begin{cor}\label{corKK}
Any two cyclic configurations of type $(v_4)$ are isomorphic if and only if these are multiplier equivalent. 
\end{cor}

\proof 
We prove the statement  for connected configurations. The general case follows then by using 
the decomposition in \eqref{decomp} and induction on $v$.  

Let $\con(\Z_v,S_r),$ $r=1,2,$ be two 
connected configurations of type $(v_4)$. Then $\langle S_r - S_r \rangle =\Z_v,$ see \eqref{conn},    
and we apply Theorem \ref{KK} to the respective  line-point incidence matrices. Now, one only needs to 
exclude the possibility that the sets $S_r$  are described by part (1b) of Theorem \ref{KK}.  That this is indeed the case 
follows from  Corollary \ref{corHMP}  where choose $H$ to be the subgroup of order $2$. \QED 

\begin{rem}
As noted in the introduction, cyclic configurations are equivalent to cyclic Haar graphs of girth $6,$ 
and therefore each of our results has a counterpart in the context of cyclic Haar graphs.  
In this spirit  \cite[Theorem 5.23]{PisS13} summarizes Corollaries \ref{corWZ} and \ref{corKK}.  
\end{rem}

\section{Proof of Theorem A}

Recall that, two configurations $\con(\Z_v,S_r),$ $r=1,2,$ are multiplier equivalent if and only if  their base lines  
$S_r$ are in the same orbit of $AGL_1(v)$. Thus Corollary \ref{corWZ} gives that $\# C(v_3)$  is equal to 
the number of orbits of $AGL_1(v)$ acting on $B_{\rm con}(v,3)$.  

\begin{lem}\label{LEM1}  
Let $v$ and $k$ be integers such that $k \ge 3$ and $v \ge k^2-k+1,$ and denote by  
$\mathcal{N}$ the  number of orbits of $AGL_1(v)$ acting on $B_{\rm con}(v,k)$. Then 
\[  
\mathcal{N} =  \frac{1}{k \phi(v)}   \sum_{l \in \Z_v^* } N(v,k,l),
\] 
where $N(v,k,l) =  \big\{ X \in B_{\rm con}(v,k) :  0 \in X \text{ and }  l  X = X - x  \text{ for some } x \in X \big\}$. 
\end{lem}

\proof  For short we put $B_0 = \{ X \in B_{\rm con}(v,k) : 0 \in X\},$ and for $X \in B_0$ with $X = \{x_1,x_2,...,x_k\},$ 
define the set 
\[ 
\widehat{X} = \{ X-x_1,X-x_2,...,X-x_k \}.
\]
It is easily seen that for every set $Y = X-x_i$ it holds $\widehat{Y} = \widehat{X}$.
It follows from this that the sets $\widehat{X}, X \in B_0,$ form a partition of $B_0$. This partition will be 
denoted by $\pi$. Notice also  that $|\widehat{X}| = k$ holds for every class $\widehat{X} \in \pi$ 
because $|X-X| = k^2 - k+1$ (see (2) in Lemma \ref{HMP}). 
Let us consider the action of $\Z_v^*$ on $B_0$ defined by $X^l = l \, X = \{ l x : x \in X\}$ for every 
$l \in \Z_v^*$ and $X \in B_0$. The partition 
$\pi$ is preserved by $\Z_v^*$ in this action, denote by  
$\orb(\Z_v^*,\pi)$ the set of the corresponding orbits.  
For $X \in B_0,$ denote by $O(X)$ the orbit of $X$ under $AGL_1(v),$ and by $O(\widehat{X})$ the orbit of $\widehat{X}$ under 
$\Z_v^*$. 

We claim that the mapping $f : O(\widehat{X}) \mapsto O(X)$  establishes  a
 bijection from $\orb(\Z_v^*,\pi)$ to the set of orbits of $AGL_1(v)$ acting on $B_{\rm con}(v,k)$ 
(notice that, the mapping $f$ is well-defined). 
It is clear that $f$ is surjective.  To settle that it is also  injective  choose $X,Y \in B_{\rm con}(v,k)$ such that $O(X) = O(Y)$. We may assume without loss of generality that  $0 \in X \cap Y$. By definition, $Y = aX+b$ for some $a \in \Z_v^*$  and $b \in \Z_v$.  Since $0 \in Y,$ $b=-ax$ for some $x \in X$. Thus $a' Y = X - x,$ where $aa' \equiv 1(\text{mod }v)$, implying that $O(\widehat{X}) = O(\widehat{Y}),$ and so $f$ is also injective, hence bijective. 
We obtain that the required number $\mathcal{N} = | \orb(\Z_v^*,\pi) |$.  Then the orbit-counting lemma applied to 
$\orb(\Z_v^*,\pi)$ yields the formula  
(see \cite[Theorem 1.7A]{DixM96}):   
\[ 
\mathcal{N}  =  \frac{1}{\phi(v)} \sum_{l \in \Z_v^*} \big| \big\{ \widehat{X} \in \pi : \widehat{X} l = \widehat{X} 
\big\} \big|.
\] 

In order to finish the proof one only needs to observe that $\widehat{X} l = \widehat{X}$ happens exactly when 
$l X = X - x$ for some $x \in X;$ and if this is so, then every set $Y \in \widehat{X}$ satisfies $lY = Y - y$ for some $y \in Y$.
This gives us 
\[
| \{ \widehat{X} \in \pi  : \widehat{X} l = \widehat{X} \} | = \frac{ N(v,k,l)}{k}.
\]
The lemma is proved. \QED 

By Corollary \ref{corWZ} and Lemma \ref{LEM1}, we find that,  
\begin{equation}\label{sum1}
\# C(v_3) = \frac{1}{3\phi(v)} \sum_{l \in \Z_v^*} N(v,3,l).
\end{equation}
We compute next the parameters $N(v,3,l)$  in \eqref{sum1}.

\medskip

Define first  the function $\Phi : \N \to \N$ by $\Phi(1)=1,$ and for $v >1$ let 
 \[ 
\Phi(v) = v\Big( 1+\frac{1}{p_1} \Big) \cdots \Big( 1+\frac{1}{p_k} \Big),
\] 
where $v$ has prime factorization $v=p_1^{n_1} \cdots p_k^{n_k}$.
Obviously, $\Phi$ is a multiplicative function, i.e., 
$\Phi(v_1v_2) =  \Phi(v_1) \Phi(v_2)$ whenever $\gcd(v_1,v_2)=1$.

\begin{lem}\label{LEM2} 
If $v > 4,$ then 
\[   
N(v,3,1) = \left\{ \begin{array}{rl} 
\frac{1}{2} \phi(v)(\Phi(v)-6)                 & \mbox{ if  $v$ is odd}, \\  [+1.5ex]

\frac{1}{2} \phi(v)(\Phi(v)-6) - 3\phi(v/2)  & \mbox{ if  $v$ is even} .
\end{array} \right. 
\] 
\end{lem}

\proof Define the sets: 
\[ 
\begin{array}{lcl}
 S(v)    &=&  \{ (x,y) \in \Z_v \times \Z_v :  \langle x,y \rangle = \Z_v \},  \\  [+1.5ex]
S^*(v) &=& \{ (x,y) \in S(v) : |\{0,x,y,-x,-y,x-y,y-x\}| < 7 \} .
\end{array}
\]

We leave for the reader to verify that the function $v \mapsto |S(v)|$ is multiplicative.
Let $v = p^n,$ $p$ is a prime. Then two elements $x,y$ generate $\Z_{v}$ if and only if one of 
them is a generator. By this we calculate  that $|S(v)| = 2\phi(v)v - \phi(v)^2 = \phi(v)(2v-\phi(v))=\phi(v) \, \Phi(v)$. 
We find, using that all functions $\phi,\Phi$ and $v \mapsto |S(v)|$ are multiplicative, that $|S(v)| = \phi(v) \, \Phi(v)$ 
for every number $v$. 

\medskip

Now, for every $x,y \in \Z_v,$ $\{0,x,y\} \in B_{\rm con}(v,3)$ if and only if $(x,y) \in S(v) \setminus S^*(v)$. Therefore, 
\begin{equation}\label{N1}
N(v,3,1) = \frac{|S(v)| - |S^*(v)|}{2} = \frac{1}{2} \big( \, \phi(v) \Phi(v) -|S^*(v)|  \, \big ). 
\end{equation}

It remains to calculate $|S^*(v)|$. 
Let $v$ be odd. Then $S^*(v)$ can be expressed as    
\[
S^*(v) = \{ (0,x), (x,0), (x,x), (x,-x), (x,2x), (2x,x) : x \in \Z_v^* \}. 
\]
Since $v > 4,$ there is no coincidence between the above pairs, and so $|S^*(v)| = 6 \phi(v)$. 
The formula for $N(v,3,1)$  follows by this and \eqref{N1}. 

Let $v$ be even, say $v=2u$. In this case 
\begin{eqnarray*}
S^*(v) & = & \{ (0,x), (x,0), (x,x), (x,-x), (x,2x), (2x,x) : x \in \Z_v^* \}  \cup \\ [+1.5ex]
            &   & \{ (u,x), (x,u), (x,x+u) : x \in \Z_v \text{ and } \langle x,u \rangle = \Z_v \}.
\end{eqnarray*}
Again, since $v > 4,$ there is no coincidence between the above pairs. 
A  quick computation gives  that $|S^*(v)| = 6\phi(v) + 6\phi(u)$.  
The formula for $N(v,3,1)$ follows by this and \eqref{N1}. The lemma is proved. \QED 

For  $l \in \Z_v^*,$ denote by $\mo(l)$ the order of $l$ as an element of $\Z_v^*$. 
Furthermore, $O(l)$ denotes the set of orbits of $\Z_v$ under $l,$  i.e., 
\[
O(l) = \big\{  \;  \{x,lx,...,l^{m-1} x\}  : x \in \Z_v \; \big\} \text{ where } m = \mo(l).
\]

\begin{lem}\label{LEM3} 
Let  $l \in \Z_v^*,$ $l \ne 1$. 
\begin{enumerate}[(i)]
\item If $\mo(l) > 3,$ then $N(v,3,l) = 0$. 
\item If $\mo(l) = 2,$ then 
\[
N(v,3,l) = \left\{ \begin{array}{cl} 
0                              & \mbox{ if } l + 1 \equiv 0(\text{mod }v),  \text{ or } v \equiv 0(\text{mod }4)  \text{ and } 
l \equiv 1(\text{mod }v/2), \\
\frac{3 \phi(v)}{2}  & \mbox{ otherwise}. 
\end{array} \right. 
\]
\item If $\mo(l)=3,$ then  
\[ 
N(v,3,l) = 
\left\{ \begin{array}{cl} 
0           & \mbox{ if } l^2+l+1 \not\equiv 0(\text{mod }v),  \\
\phi(v)  & \mbox{ otherwise}. 
\end{array} \right. 
\]
\end{enumerate}
\end{lem}

\proof  Put again $B_0 = \{X \in B_{\rm con}(v,3) : 0 \in X\},$ and let $X \in B_0$ such that $X = \{0,x,y\}$ and 
\begin{equation}\label{NewE1}
l X = X \text{ or }  l X = X -x.
\end{equation}
   
We consider step-by-step all cases (i)-(iii). 

\medskip

(i):  Assume by contradiction that \eqref{NewE1} holds for some $l \in \Z_v^*$ with $\mo(l) > 3$. 
If $X l = X,$ then $l^2 x = x$ and $l^2 y = y$.  This together with $\langle x,y \rangle =\Z_v$ imply that 
$l^2 \equiv 1(\text{mod }v),$  a contradiction to  $\mo(l) > 2$. 
Let $X l = X - x,$ and so $\{l x,l y\} = \{-x,y-x\}$. Now, if $l x = -x$ and $l y =  y - x,$ then 
$l^2 x =x$ and $l^2 y = y$ which is impossible. If $l x = y-x$ and $l y = -x,$ then $l^3 x = x$ and $l^3 y = y,$ 
implying that $l^3 \equiv 1(\text{mod }v),$ which is  in contradiction with  $\mo(l) > 3$. 

\medskip

(ii):  Assume that \eqref{NewE1} holds with $\mo(l) = 2$. 
If $l X = X,$ then $l x = y$ and $l y = x$ and so we find $X$ as 
$X = \{0,x,lx\} , \;  x \in \Z_v^* $. Let $lX = X - x$. Then it follows that $l x = -x$ and $l y = y-x$ 
(otherwise $l^3 \equiv 1(\text{mod }v),$ a contradiction to $\mo(l) = 2$), and so 
$X = \{0,y,-ly+y\}$ where $y \in \Z_v^*$.
Since $X \in B_0,$ the elements $0, 1, -1, l, -l, l-1 \text{ and } 1-l $ must be pairwise distinct. 
We conclude from these that, $N(v,3,l) = 0$ if $l + 1 =  0(\text{mod }v)$ or   
$l  \equiv  1(\text{mod }v/2),$ and otherwise $N(v,3,l) $ is the size of the following set:
\[
 \big\{  \{ 0,x,lx \} : x \in \Z_v^*  \big\}   \; \cup \; \big\{  \{ 0,x,-lx+x \} : x \in \Z_v^*  \big\}.
\]
We observe in turn that, the two sets above are disjoint, the first has size $\phi(v)/2,$ 
while the second has cardinality $\phi(v)$. Then  (ii) follows.  

\medskip

(iii): Assume that \eqref{NewE1} holds with $\mo(l) = 3$. Then 
$X = lX - x,$  $lx = y-x$ and $ly = -x$ (otherwise $l^2 \equiv 1(\text{mod }v),$ see above). 
Thus $X = \{0,x,x+lx\},$ $x \in \Z_v^*$ and $l^2+l \equiv -1(\text{mod }v)$. We conclude that, 
$N(v,3,l) = 0$ if $l^2 + l + 1 \not\equiv  0(\text{mod }v),$ and  otherwise 
$N(v,3,l)  =  \big| \big\{  \{ 0,x,lx+x \} : x \in \Z_v^*  \big\} \big| = \phi(v)$.
Thus (iii) follows, and this completes the proof of the lemma. \QED 

\medskip

\noindent{\sc Proof of Theorem A:} By Lemmas \ref{LEM1} and \ref{LEM2}, the sum in \eqref{sum1} reduces to 
\begin{equation}\label{sum2}
\# C(v_3) = \left\{ \begin{array}{lr}
\frac{1}{6} \Phi(v)  - 1 + \frac{1}{2} \gamma_1 + \frac{1}{3} \gamma_2  & \text{ if $v$ is odd,} \\  [+1.5ex]
\frac{1}{6} \Phi(v)  - \frac{\phi(v/2)}{\phi(v)}  - 1  + \frac{1}{2} \gamma_1 + \frac{1}{3} \gamma_2  & \text{ if $v$ is even,}
                      \end{array} \right.
\end{equation}
where $\gamma_1$ and $\gamma_2$ are defined by 
\[ 
\begin{array}{lcl}
\gamma_1    & = &  | \{ l \in \Z_v^* : \mo(l) = 2, \;  l+1 \not\equiv 0(\text{mod }v)  \text{ and } 
                                l \not\equiv 1(\text{mod }v/2) \text{ if } v \equiv 0(\text{mod }4) \} |,  \\  [+1.5ex]
\gamma_2    & = & | \{ l \in \Z_v^* : \mo(l) = 3 \text{ and } l^2+l+1 \equiv 0(\text{mod }v)\} |.
\end{array}
\]

In calculating  $\gamma_1$ and $\gamma_2$ below we shall use the fact  $\Z_v^*$ 
can be written as 
$\Z_v^*  = \Z_{p_1^{n_1}}^* \times \cdots \times \Z_{p_k^{n_k}}^*,$ and  every $l \in \Z_v^*$ can be expressed as 
\begin{equation}\label{li}
l = ( l_1, ..., l_k), \text{ where } l_i \in \Z_{p_i^{n_i}}^*  \text{ for  every } i \in \{1,...,k\}.
\end{equation}
Note that, we may assume that $l_i \equiv l(\text{mod }p_i^{n_i})$ for every $i \in \{1,...,k\}$. 

\bigskip

\noindent {\sc Case 1.} $v$ is odd.

\medskip

Since $v$ is odd, there are exactly $2^k -1$ elements $l \in \Z_v^*$ such that $\mo(l) = 2,$ and all but one 
contributes to $\gamma_1$ (namely, $l = v-1$ is excluded in the definition of $\gamma_1$). 
Thus $\gamma_1 = 2^k - 2$. The value of $\gamma_2$ depends solely on the residue of 
$v$ modulo $9$ and  the reside of prime factors $p_i$ modulo $3$. 
Let $l \in \Z_v^*$ such that $\mo(l) = 3$ and write $l = (l_1,...,l_k)$ as described in \eqref{li}.
Thus $l_i$ is of order $1$ or $3$ in $\Z_{p_i^{n_i}}^*$. 

\medskip

{\sc Case 1.1.} $p_i \equiv 1(\text{mod }3)$ for every $i \in \{1,...,k\}$ . 

\medskip

If $l_i$ is of order $1$ in $\Z_{p_i^{n_i}}^*,$ then  
$l \equiv l_i \equiv 1(\text{mod }p_i^{n_i}),$ from which $l^2 + l + 1 \equiv 3(\text{mod }p_i^{n_i}),$ 
hence $l^2 + l + 1 \not\equiv 0(\text{mod }v),$ so $l$ cannot contribute to $\gamma_2$. If $l_i$ is of order $3$ in $\Z_{p_i^{n_i}}^*,$ then $l^2 + l + 1 \equiv l_i^2 + l_i +1 \equiv 0(\text{mod }p_i^{n_i})$ for every $i \in \{1,...,k\},$ hence 
$l^2 + l + 1 \equiv 0(\text{mod }v)$. Since there are exactly two elements in $\Z_{p_i^{n_i}}$ of order $3,$ 
$\gamma_2 = 2^k$. Substitute this and $\gamma_1 = 2^k-2$ in \eqref{sum2}. We obtain that 
$\# C(v_3) = \frac{1}{6} \Phi(v)  + \frac{5}{6} 2^k - 2$.  

\medskip

{\sc Case 1.2.}  $v \equiv 3(\text{mod }9)$ and $p_i \equiv 0/1(\text{mod }3)$ for every $i \in \{1,...,k\}$ . 

\medskip

We may write $p_1^{n_1} = 3$. We obtain, by the same argument as  in the previous case, that 
$l$ contributes to $\gamma_2$ if and only if $l_1$ is of order $1$ in $\Z_{p_1^{n_1}},$ and $l_i$ is of order $3$ in $\Z_{p_i^{n_i}}^*$ if $i \ge 2$. Thus $\gamma_2 = 2^{k-1},$ which together with $\gamma_1 = 2^k-2$ yield  in 
\eqref{sum2} that $\# C(v_3) = \frac{1}{6} \Phi(v)  + \frac{2}{3} 2^k - 2$. 

\medskip

{\sc Case 1.3.}  $v \equiv 0(\text{mod }9)$ or $p_i \equiv 2(\text{mod }3)$ for some $i \in \{1,...,k\}$ . 

\medskip

We show that in this case $l^2 + l +1 \not\equiv 0(\text{mod }v)$ independently of the choice $l$.  
Thus $\gamma_2  = 0,$ and so $\# C(v_3) = \frac{1}{6} \Phi(v)  + \frac{1}{2} 2^k - 2$. 

Suppose first that $v \equiv 0(\text{mod }9)$.  We may write $p_1 = 3,$ now  $n_1 \ge 2$. 
Since $\mo(l)=3,$ $l_1 \equiv 1(\text{mod } 3^{n_1-1})$. We claim that $l_1^2 + l_1 + 1 \equiv 3(\text{mod } 3^{n_1})$. 
Indeed, $l_1 \equiv 3^{n_1-1}k+1 (\text{mod }3^{n_1})$ for some $k \in \{0,1,2\}$. Hence 
\[
l_1^2 + l_1 + 1 \equiv   (k+2k)3^{n_1-1} + 3 \equiv 3(\text{mod }3^{n_1}).
\]
Therefore, $l^2+l+1 \equiv l_1^2 + l_1 + 1 \equiv 3(\text{mod } 3^{n_1}),$ and since 
$n_1 \ge 2,$ $l^2+l+1 \not\equiv 0(\text{mod }3_1^{n_1}),$ and so $l^2+l+1 \not\equiv 0(\text{mod }v)$.

Suppose next that $p_i \equiv 2(\text{mod }3)$ for some $i \in \{1,...,k\}$. Then $l_i$ must be of order $1$ in $\Z_{p_i^{n_i}},$ 
and hence $l^2+l+1 \equiv l_i^2 + l_i + 1 \equiv 3(\text{mod }p_i^{n_i}),$ and so 
$l^2+l+1 \not\equiv 0(\text{mod }v)$. 
 
\bigskip

\noindent {\sc Case 2.} $v$ is even.

\medskip

Since $v$ is even, $l$ is odd, and thus $l^2+l+1 \not\equiv 0(\text{mod } v)$. 
We obtain that $\gamma_2 = 0$. The value of $\gamma_1$ depends on the residue of $n$ modulo 
$8$. The number of elements of order $2$ in $\Z_v^*$ is $2^{k-1}-1$ if $v \equiv 2/6(\text{mod } 8),$ 
$2^k-1$ if $v \equiv 4(\text{mod }8),$ and $2^{k+1}-1$ if $v \equiv 0(\text{mod }8)$ (see \cite[Exercise 6.12]{JonJ98}). Thus  
\begin{equation}\label{gamma1}
\gamma_1 = \left\{ \begin{array}{rcl}
2^{k-1} - 2    & \text{ if } & v \equiv 2/6(\text{mod }8), \\ [+1.5ex] 
2^k - 3          & \text{ if } & v \equiv 4(\text{mod }8),  \\ [+1.5ex]
2^{k+1} - 3  & \text{ if } & v \equiv 0(\text{mod }8). 
\end{array} \right.
\end{equation}

Obviously, $\phi(v/2)/\phi(v) = 1$ if $v \equiv 2(\text{mod }4)$ and it is $1/2$ if $v \equiv 0(\text{mod }4)$.
Substituting this, \eqref{gamma1} and $\gamma_2 = 0$ in  \eqref{sum2} yields formula 
\eqref{formula}. The theorem is proved. \QED 

\section{Proof of Theorem B}

We consider cyclic configurations in the wider context of cyclic  objects, where  
by a \textit{cyclic object} of order $v$ we mean a relational structure on $\Z_v$ 
which is  invariant under the translation $\tau : x \mapsto x+1$. 
The  set of all cyclic objects of order $v$ will be denoted by $\obj(\tau,\Z_v)$  (see \cite{Muz99}). 
An \textit{isomorphism} between two cyclic objects $X_r,  r = 1,2,$ is a permutation  
$\sigma$ of $\Z_v$ which maps $X_1$ onto $X_2,$ if $X = X_1  = X_2,$ then 
$\sigma$ is an \textit{automorphism}, the group of all automorphisms will be denoted by $\aut(X)$.   
Given a class $\K$ of objects in $\obj(\tau,\Z_v),$ a \textit{solving set} for $\K$  is a  set $\Delta$ of  permutations of 
$\Z_v$ satisfying the following property (see \cite{Muz99}): 
\[
(\forall X \in \K) \; (\forall Y \in \obj(\tau,\Z_v))  \;  (X \cong Y  \iff   X^\sigma = Y \text{ for some } \sigma \in \Delta). 
\]

P\'alfy's  characterization of CI-groups (see the paragraph before Theorem B) yields the following theorem:

\begin{thm}\label{P}
{\rm \cite{Pal87}} 
The set $\Z_v^*$ is  a solving set for $\obj(\tau,\Z_v)$ if and only if $v = 4$ or  
$\gcd(v,\phi(v)) = 1$.\footnote{Here $\Z_v^*$ denotes the set of permutations $x \mapsto ax$ where $a$ goes over the set of all units in $\Z_v$.}   
\end{thm}

Let  $p$ and $q$ be distinct primes. For every object $X \in \obj(\tau,\Z_{pq}),$ a solving set for $X$ was determined by Huffman \cite{Huf96}.  Before stating the relevant results, let us recall the required notations.
For $j \in \Z_v^*,$ let $\mu_j$ be the permutation 
$\mu_j : x \mapsto jx$. For $i \in \{0,1,...,q-1\},$ define the permutation $\tau_i$ by 
\[
\tau_i:       x  \mapsto \begin{cases} x+q & \text{if } x \equiv i(\text{mod }q), \\ x & \text{otherwise},  \end{cases} 
\]
and if in addition $j \in \Z_v^*$ with $j \equiv 1(\text{mod }q),$ then define the permutation $\mu_{i,j}$ by    
 \[
\mu_{i,j} : x  \mapsto \begin{cases} jx & \text{if } x \equiv i(\text{mod }q), \\ x & \text{otherwise}. \end{cases}
\] 

For the next two theorems suppose in addition that $q$ divides $p-1$. Furthermore,  fix an element $a \in \Z_v^*$ of 
order $\mo(a) = p-1$ for which $a \equiv 1(\text{mod }q),$ and put 
$b = a^{(p-1)/q}$.

\begin{thm}\label{H1} 
{\rm \cite[Theorem 1.1]{Huf96}} Let $v = pq,$ where $p,q$ are primes such that $q$ divides $p-1,$ and let 
$X \in \obj(\tau,\Z_v)$ such that $\mu_b  \notin \aut(X),$ where $b$ is defined above. 
Then $\Z_v^*$ is a solving set for $X$.
\end{thm}

The powers $a,a^2,..,a^p$ are pairwise distinct modulo $p$. Let $\alpha$ be the positive integer in $\{1,2,...,p\}$ 
that $a^\alpha \equiv -s(\text{mod }p),$ where $s = (p-1)/q$.  For $i \in \{0,1,...,q-1\}$ define 
$\nu_i = \prod_{j = 0}^{q-1} \mu_{j,a^\alpha b^{-ij}}$.
Notice that, $\nu_0 = \mu_{a^\alpha}$.
The next theorem is \cite[Theorem 1.2]{Huf96}, which, for our convenience, is formulated  slightly 
differently. 

\begin{thm}\label{H2} 
Let $v = pq,$ where $p,q$ are primes such that $q$ divides $p-1,$ and let 
$X \in \obj(\tau,\Z_v)$ such that $\mu_b  \in \aut(X)$ and $\tau_0 \notin \aut(X),$ where $b$ is defined above. 
Let $\beta$ be the smallest positive integer such that $\mu_a^\beta \in \aut(X)$.
Then $X$ admits a  solving set $\Delta$ in the form:   
\begin{equation}\label{Delta}
\Delta = \Big\{ \mu_a^i \, \nu_k \, \mu_j^{-1} :  0 \le i <  \beta, \; 0 < j  \le q-1, \; 0 \le k \le q-1, \; 
\prod_{l=0}^{q-1}\tau_l^{b^{(l+1)k}} \in  \aut(X)  \Big\}.
\end{equation}
\end{thm}

The last result before we prove Theorem B is a special case of \cite[Lemma 3.1]{Bab77}.  

\begin{lem}\label{B}
The following (1)-(2) are equivalent for every object $X \in \obj(\tau,\Z_v)$.
\begin{enumerate}[(1)]
\item $\Z_v^*$ is a solving set for $X$.
\item Every two regular cyclic subgroup of $\aut(X)$ are conjugate in $\aut(X)$.
\end{enumerate}
\end{lem}

\bigskip

\noindent{\sc Proof of Theorem B:} \ Obviously, the theorem can be rephrased as follows: 
$\Z_v^*$ is a solving set for the class of  cyclic configurations on $v$ points if $v = pq$ or $v = p^n,$ where $p,q$ are primes.

\medskip

\noindent{\sc The case $v = pq$:} We prove the above statement  for connected configurations. The general case follows then by using the decomposition in \eqref{decomp} and the fact that the statement is true for configurations with a  prime number of 
points, so let $\C = \con(\Z_{pq},S)$ be a connected cyclic configuration. 

Towards a contradiction assume that $\Z_{pq}^*$ is not a solving set for $\C$. 
Because of Theorem \ref{P} we may also assume that $q$ divides $p-1$. In the rest of the proof we keep the previous notations: $\tau_0,a,b,\alpha$  and $\nu_0,\nu_1,...,\nu_{q-1}$. Let 
$P = \{0,q,...,(p-1)q\},$ i.e., the subgroup of $\Z_{pq}$ of order $p$. 
Replace $S$ with a suitable line $S+i$ if necessary to ensure that $S \cap P \ne \emptyset$.  Also, 
$S \not\subseteq P$ by the connectedness of $X,$ i.e., there exists $t \in \{1,...,q-1\}$ such that 
\begin{equation}\label{nonempty}
S \cap P \ne \emptyset \text{ and } S \cap (P+ t) \ne \emptyset.
\end{equation}

Suppose for the moment that $\tau_0 \in \aut(\C)$. Using that $\tau_0$ fixes every point outside $P,$ \eqref{nonempty} 
and that $|S| \ge 3,$ we conclude $|S^{\tau_0^k} \cap S| \ge 2$ for some $k \in \{1,...,q-1\}$.  Hence $S^{\tau_0^k} = S$. 
As $P$ is an orbit of $\tau_0^k,$ $P \subseteq S,$ which contradicts Corollary \ref{corHMP}.
Thus $\tau_0 \notin \aut(\C)$. 

Therefore, Theorems \ref{H1} and \ref{H2}, together with the assumption that $\Z_{pq}^*$ is not a solving set,  
imply that $\C$ admits a solving set $\Delta$ defined in \eqref{Delta}. 
Consider the permutation $\sigma = \prod_{l=0}^{q-1}\tau_l^{b^{(l+1)k}},$ $k \in \{0,1,...,q-1\}$. 
If $k = 0,$ then $\sigma = \tau^q$ which is clearly in $\aut(\C)$. The 
corresponding permutations in $\Delta$ are  $\mu_a^i \, \nu_0 \, \mu_j^{-1} = 
\mu_a^i \, \mu_{a^\alpha} \, \mu_j^{-1}$. 
Since $\Delta \not\subseteq \Z_v^*,$ there must exist $k > 0$ for which
$\sigma = \prod_{l=0}^{q-1}\tau_l^{b^{(l+1)k}}$ belongs to $\aut(\C)$. 
Notice that, 
\begin{equation}\label{distinct}
\forall i,j \in \{0,1,...,q-1\}:  \; i \ne j \implies b^{ik} \not\equiv b^{jk} (\text{mod }p).
\end{equation}
For otherwise, $b^{(i-j)k} \equiv 1(\text{mod }p)$. Since $\mo(a) = p-1,$ $a \equiv 1(\text{mod }q)$ and 
$b =  a^{(p-1)/q},$ we find from $a^{(p-1)(i-j)k/q} = b^{(i-j)k} \equiv 1(\text{mod }p)$ that 
$p-1$ divides  $(p-1)(i-j)k/q,$ and so $q$ divides $(i-j)k,$  a contradiction. 

Consider the product $\sigma' = \sigma \tau^{-b^k}$. Now, $\sigma'$ fixes each point in $P,$ but because of \eqref{distinct} 
it permutes the points of $P+t$ in a $p$-cycle. 
Unless $|S \cap (P+x)| \le 1$ for every $x \in \{0,1,...,q-1\},$ we may also assume that 
$|S \cap P| \ge 2$. However, if $|S \cap P| \ge 2,$ then $\sigma'$ fixes $S,$ implying that $(P+t) \subseteq S,$ which is 
impossible. 

We are left with the case that $|S \cap (P+x)| \le 1$ for every $x \in \{0,1,...,q-1\}$. 
Note that, then the same holds for all lines $S+i$. It is obvious that $|S| \le q$.  Let $\{s\} =  S \cap P$. As $\C$ is balanced, there are exactly $|S|$ lines through $s$. Now, each of the  lines $S,S^{\sigma'},...,S^{\sigma'^{p-1}}$ contains $s$, while they intersect $P+t$ at distinct points. These imply in turn that, they are pairwise distinct, hence $|S| \ge p$, and so 
$p \le  |S| \le q,$ a contradiction.    This completes the proof of case $v=pq$.

\medskip

We turn next to the case  $v=p^n$. Now, we cannot relay on a list of solving sets covering all cyclic objects as such 
list is available only when $v = p^2$ (see  \cite{HufJP93}). The argument below will be a combination of 
Lemma \ref{B} with Sylow's theorems.

\bigskip

\noindent {\sc The case $v = p^n$:} Again, it is sufficient to consider connected configurations, the general case follows then by using the decomposition in \eqref{decomp} and induction on $n$.  
Let $\C = \con(\Z_{p^n},S)$ be a connected cyclic configuration, $G=\aut(\C)$ and $C$ be the group generated by 
$\tau : x \mapsto x+1$. Let $G_p$ be a Sylow $p$-subgroup of $G$ such that $C \leq G_p$. By Lemma \ref{B} and 
Sylow's theorems it is sufficient to prove that $G_p = C$.

Towards a  contradiction assume that $C < G_p$. Then the normalizer $N_{G_{p}}(C) >  C$. Let us put $N = N_{G_{p}}(C)$ and let $N_{0}$ be the stabilizer of $0$ in $N$. Then $N_{0}$ is non-trivial, and we may choose  $\sigma$ from $N_{0}$ of order $p$. 
Since $\sigma$ normalizes the regular subgroup $C$ and fixes $0$, $\sigma = \mu_a$ for some $a \in \Z_{p^n}^*$ 
(see \cite[Exercise 2.5.6]{DixM96}).  Then $\mo(a) = p$. 
Using the well-known structure of $\Z_{p^n}^*$ (cf. \cite[Theorem 6.7 and Exercise 6.12]{JonJ98}) we deduce that  $n \ge 2,$ and either   
\[ 
a = a'p^{n-1}+1 \text{ for some  } a' \in \{ 1, \dots , p-1 \}, 
\] 
or $n \ge 3,$ $p = 2$ and $a \in \{ 2^n -1,2^{n-1} - 1 \}$.

Assume for the moment that the latter case holds. 
Let $Q = \langle C,\, \sigma \rangle$. It is a routine exercise to show that $C$ is the only cyclic subgroup of $Q$ of order $2^n$. 
This implies that the normalizer $N_{G_{2}}(Q) \leq N_{G_2}(C) = N$. Let $H$ be an arbitrary regular cyclic subgroup of $G$. If $Q = G_2$, then, by Sylow's theorems, $H^g < Q$ for some $g \in G$, and so $H^g = C,$ and we are done by Lemma \ref{B}. 
Thus we may assume that $Q < G_2$.
Then $Q <  N_{G_{2}}(Q) \leq N$. Choose an element $\sigma' \in N_0$ such that $\sigma' \ne \sigma$. 
It is well-known that $5^{2^{n-3}} \equiv 2^{n-1}+1(\text{mod }2^n)$ (see \cite[Lemma 6.9]{JonJ98}), 
and that $\Z_{2^n}^* = \langle 5 \rangle \times \langle - 1 \rangle \cong \Z_{2^{n-2}} \times \Z_2$ 
(see \cite[Theorem 6.10]{JonJ98}).  These imply that 
$\mu_{2^{n-1}+1} \in \langle \sigma,\sigma' \rangle,$ and so $\mu_{2^{n-1}+1} \in N_0$.  
Therefore, we may assume that $\mu_{a} \in \aut(\C)$ where $a = a'p^{n-1}+1$ for some $a' \in \{ 1, \dots , p-1 \}$.

Now, $\mu_a$ maps $S$ to a line of $\C,$ hence we may write $a S + b = S$ for some $b \in \Z_{p^n}$.  
Equivalently, $S$ is a union of orbits of the affine transformation $\varphi : x \mapsto ax + b$. 
Then $\varphi^p$ is equal to the translation $x \mapsto x + (1 + a + \cdots + a^{p-1})b$. 
By Corollary \ref{corHMP}, $S$ contains no non-trivial cosets. Form this and that $S$  is a union orbits of $\varphi^p$, 
we find  that $(1 + a + \cdots + a^{p-1})b  \equiv 0 (\text{mod } p^n)$. This quickly implies that $p^{n-1}$ divides $b,$ 
hence we may write $b = b'p^{n-1}$ for some $b' \in \{0,1,...,p-1\}$.  Also, 
 \[
\varphi : x \mapsto ax + b = x + (a'x + b')p^{n-1}.
\]
From this we easily find the orbits of $\varphi$. For $x \in \Z_v,$ let $O$ be the orbit which contains  $x$. Then 
\[
O =  \begin{cases} \{x\}  & \text{ if } a'x + b' \equiv0(\text{mod }p), \\ 
                               P + x   & \text{otherwise},
        \end{cases}
\]
where $P = \{0, p^{n-1},...,(p-1)p^{n-1}\},$ i.e., the subgroup of $\Z_{p^n}$ of order $p$. 
Since $X$ is connected, $\langle S - S \rangle = \Z_{p^n}$. This implies 
that $a's + b' \not\equiv 0(\text{mod } p )$ for some $s \in S$. 
But then the coset $(P + s) \subseteq S,$ contradicting Corollary \ref{corHMP}. The theorem is proved  \QED

\end{document}